\documentstyle[12pt]{article}
\title{Periodical Solutions of \\Multi-Time Hamilton Equations}
\author{Constantin Udri\c ste, Iulian Duca}
\date{}
\begin{document}
\maketitle

\newcommand{\di}{\displaystyle}
\newcommand{\noa}{\noalign{\medskip}}

\begin{abstract}

To our knowledge, there are two main references [9], [12]
regarding the periodical solutions of multi-time Euler-Lagrange
systems, even if the multi-time equations appeared in 1935, being
introduced by de Donder. That is why, the central objective of
this paper is to solve an open problem raised in [12]: what we can
say about periodical solutions of multi-time Hamilton systems when
the Hamiltonian is convex?

Section 1 recall well-known facts regarding the equivalence
between Euler-Lagrange equations and Hamilton equations. Section 2
analyzes the action that produces multi-time Hamilton equations,
and introduces the Legendre transform of a Hamiltonian together a
new dual action. Section 3 proves the existence of periodical
solutions of multi-time Hamilton equations via periodical
extremals of the dual action, when the Hamiltonian is convex.
\end{abstract}

\textbf{Key words:} multi-time Hamilton action, periodical extremals, convex Hamiltonian.

\textbf{2000 Mathematics Subject Classification:} 53C21, 49K20, 49S05.

\section{Classical Legendre transform of a Lagrangian}

In this Section we recall the classical duality used in mechanics, based on the fact that the Legendre transform of a Lagrangian is a Hamiltonian
(see [2], [7], [8], [11] for single-time theory, and [1], [3]-[6], [8]-[12] for multi-time theory).

\subsection{Single-time Euler-Lagrange equations}

The Lagrangian, $L:R\times R^{2n}\rightarrow R$, $\left( t,x,\dot x
\right) \rightarrow L\left( t,x,\dot x\right) $, where
$x=\left(x^{1}\left( t\right) ,...,x^{n}\left( t\right) \right) $ and
$\dot x =\left( \dot x^{1}\left( t\right) ,...,\dot x^{n}\left( t  \right)\right)$,
produces the single-time Euler-Lagrange equations
$$
\frac{d}{dt}\frac{\partial L}{\partial \dot x^{i}}=
\frac{\partial L}{\partial x^{i}}, \quad
i=1,...,n. \eqno (1)
$$

\subsection{Single-time Hamilton equations}

We suppose that the Lagrangian $L$ defines the diffeomorphism $\dot x^{i}\rightarrow p_{i}=\di\frac{\partial L}{\partial \dot x^{i}}$
called the {\it Legendre transformation}.

The Legendre transformation of the Lagrangian $L$ is the Hamiltonian
$H:R\times R^{2n}\rightarrow R$, $\left( t,x,p\right) \rightarrow H\left(t,x,p\right)$, $H\left( t,x,p\right) =\dot x^{i}p_{i}-L\left( t,x,\dot x\right)$.

The equations (1) become
$$
\frac{dx^{i}}{dt}=\frac{\partial H}{\partial p_{i}},\quad \frac{dp_{i}}{dt}=-\frac{\partial H}{\partial x^{i}}, \quad i=1,...,n \eqno(2)
$$
because
$$
\frac{\partial H}{\partial x^{i}}=p_{j}\frac{\partial \dot x
^{j}}{\partial x^{i}}-\frac{\partial L}{\partial x^{i}}-\frac{\partial L}{\partial \dot x^{j}}\frac{\partial \dot x^{j}}{\partial x^{i}}=
-\frac{\partial L}{\partial x^{i}}=-\frac{d}{dt}p_{i}
$$
and
$$
\frac{\partial H}{\partial
p_{i}}=\dot x^{i}+p_{j}\frac{\partial \dot x^{j}}{\partial
p_{i}}-\frac{\partial L}{\partial \dot x^{j}}\frac{\partial
\dot{x}^{j}}{\partial p_{i}}=\dot x^{i}. \eqno(3)
$$
So, the second order ODEs (1) in an $n$-dimensional space are equivalent with the first order ODEs (2) in a $2n$-dimensional space.

The Hamiltonian $H$ associated to $L$ was built in 1834 by Hamilton. With this function, the Euler-Lagrange equations (1) achieve a
symmetrical structure in the form of the Hamilton equations (2).

More than that, in the case when the Lagrangian $L$ is an autonomous function (does not explicitly depend on $t$), the function $H$ is a
first integral for the Hamilton equations, because
$$
\frac{dH}{dt}=\frac{d}{dt}\left( p_{i}
\dot x^{i}-L\right) =\frac{dp_{i}}{dt}\dot x^{i}+p_{i}\frac{d\dot x^{i}}{dt}-\frac{\partial L}{\partial t}-\frac{\partial L}{\partial x^{i}}\frac{dx^{i}}{dt}-\frac{\partial L}{\partial \dot x^{i}}\frac{d\dot x^{i}}{dt}=-\frac{\partial L}{\partial t}=0.
$$

Using the symplectic structure $J=\left(
\begin{array}{cc}
0 & \delta_{j}^{i} \\ \noa
-\delta_{j}^{i} & 0
\end{array} \right)$, $\left( \delta_{j}^{i}\right) \in M_{n}\left( R\right)$, we may write the Hamilton equations in the matrix form
$$
J\left( \begin{array}{c}
\dot x^{j} \\ \noa
\dot{p}_{i} \end{array} \right) +\left(
\begin{array}{c}
\di\frac{\partial H}{\partial x^{i}} \\ \noa
\di\frac{\partial H}{\partial p_{i}}
\end{array} \right) =\left(
\begin{array}{c} 0 \\ \noa  0 \end{array}\right).
$$

\subsection{Multi-time Euler-Lagrange equations}

We consider the multi-time variable $t=\left( t^{1},...,t^{p}\right) \in R^{p}$, the functions
$x^{i}:R^{p}\rightarrow R$ , $\left( t^{1},...,t^{p}\right)\rightarrow x^{i}\left( t^{1},...,t^{p}\right) $, $i=1,...,n$ , and we
denote $x_{\alpha }^{i}=\di\frac{\partial x^{i}}{\partial t^{\alpha }}$ , $\alpha =1,...,p$. The Lagrange function
$$
L:R^{p+n+np}\rightarrow R,  \quad
\left( t^{\alpha },x^{i},x_{\alpha }^{i}\right) \rightarrow L\left(
t^{\alpha },x^{i},x_{\alpha }^{i}\right)
$$
gives the Euler-Lagrange equations
$$
\frac{\partial}{\partial t^{\alpha }}
\frac{\partial L}{\partial x^i_{\alpha}}=\frac{\partial L}{\partial x^{i}}, \quad i=1,...,n, \quad \alpha =1,...,p
$$
(second order PDEs system on the n-dimensional space).

\subsection{Multi-time Hamilton equations}

In 1935 de Donder [1] obtained the Hamilton equations in the
multi-time case by using the partial derivatives
$$
p_{k}^{\alpha }=\frac{\partial L}{\partial x_{\alpha }^{k}}
\eqno (4)
$$
and the Hamiltonian $H=p_{k}^{\alpha }x_{\alpha }^{k}-L$. If \ $L$ satisfies some conditions, then the system (4)
defines a $C^{1}$ bijective transformation $x_{\alpha }^{i}\rightarrow p^{\alpha}_{i}$ , called the Legendre transformation
for the multi-time case. By this transformation we have:
$$
\frac{\partial H}{\partial p_{i}^{\alpha }}=x_{\alpha }^{i}+p_{k}^{\beta }\frac{\partial x_{\beta }^{k}
}{\partial p_{i}^{\alpha }}-\frac{\partial L}{\partial x_{\beta }^{k}}\frac{\partial x_{\beta }^{k}}{\partial p_{i}^{\alpha }}=x_{\alpha }^{i}
$$
$$
\frac{\partial H}{\partial x^{i}}
=p_{k}^{\alpha }\frac{\partial x_{\alpha }^{k}}{\partial x^{i}}-\frac{\partial L}{\partial x^{i}}-\frac{\partial L}{\partial x_{\alpha }^{k}}
\frac{\partial x_{\alpha }^{k}}{\partial x^{i}}=-\frac{\partial L}{\partial x^{i}}.
$$
The $np + n$ Hamilton equations $\di\frac{\partial x^{i}}{\partial t^{\alpha }}=\di\frac{\partial H}{\partial p_{i}^{\alpha }}$,
$\di\frac{\partial p_{i}^{\alpha }}{\partial t^{\alpha }}=-\di\frac{\partial H}{\partial x^{i}}$ (summation after
$\alpha $), $i=1,...,n$, $\alpha =1,...,p$ are first order PDEs on the space $R^{n+pn}$,
equivalent to the Euler-Lagrange equations on $R^{n}$.

\subsection{The conservation of energy-moment tensor}

The multi-time Hamiltonian is not conserved on the solutions of the
multi-time Hamilton equations even if the Lagrangian is autonomous (does not depend on \ $t^{\alpha }$). But, we may observe that the
 Lagrangian $L$ defines the energy-moment tensor
$$
T_{\beta }^{\alpha }=x_{\beta}^{i}\frac{\partial L}{\partial x_{\alpha }^{i}}-\delta _{\beta }^{\alpha }L,
$$
whose divergence is
$$
\frac{\partial }{\partial t^{\alpha }}T_{\beta
}^{\alpha }=\frac{\partial }{\partial t^{\alpha }}\left( x_{\beta
}^{i}p_{i}^{\alpha }\right) -\delta _{\beta }^{\alpha }\left( \frac{\partial
L}{\partial t^{\alpha }}+\frac{\partial L}{\partial x^{i}}x_{\alpha }^{i}+\frac{\partial L}{\partial x_{\gamma }^{i}}
\frac{\partial x_{\gamma }^{i}}{\partial t^{\alpha }}\right) =
$$
$$
=\frac{\partial x_{\beta }^{i}}{\partial
t^{\alpha }}p_{i}^{\alpha }+x_{\beta }^{i}\frac{\partial p_{i}^{\alpha }}{\partial t^{\alpha }}-\frac{\partial L}{\partial t^{\beta }}-
\frac{\partial L}{\partial x^{i}}x_{\beta }^{i}-\frac{\partial L}{\partial x_{\gamma }^{i}}
\frac{\partial x_{\gamma }^{i}}{\partial t^{\beta }}=
$$
$$
=x_{\beta }^{i}\left( \frac{\partial p_{i}^{\alpha }}{\partial t^{\alpha }}-\frac{\partial L}{\partial x^{i}}\right) -
\frac{\partial L}{\partial t^{\beta }}=-\frac{\partial L}{\partial t^{\beta }}.
$$
So, the energy-moment
tensor $T_{\beta }^{\alpha }$ is conserved on the solutions of the
multi-time Hamilton equations, if the Lagrangian $L$ is autonomous.

\section{Legendre transform of a multi-time Hamiltonian. Dual action}

We may write the multi-time Hamilton equations in the form ([9],[12])
$$
\delta_{\beta }^{\alpha }\delta _{j}^{i}\frac{\partial p_{i}^{\beta }}{\partial t^{\alpha }}+\frac{\partial H}{\partial x^{j}}=0, \quad -
\delta_{\beta }^{\alpha}\delta _{j}^{i}\frac{\partial x^{j}}{\partial t^{\alpha }}+\frac{\partial H}{\partial p_{i}^{\beta }}=0
$$
or
$$
\left( \delta \otimes J\right) \left(
\begin{array}{c}
\di\frac{\partial x^j}{\partial t^{\alpha }}
 \\ \noa
\di\frac{\partial p_{i}^{\beta }}{\partial t^{\alpha }}
\end{array}
\right) + \left(
\begin{array}{c}
\di\frac{\partial H}{\partial x^{j}} \\ \noa
\di\frac{\partial H}{\partial p_{j}^{\beta }}
\end{array}
\right)=
\left(
\begin{array}{c}
0 \\ \noa
0\end{array}
\right), \eqno(5)
$$
$j=1,...,n$; $i=1,...,n$; $\alpha =1,...,p$; $\beta =1,...,p$,
where
$$
\delta \otimes J=\left(
\begin{array}{cc}
0 & \delta _{\beta }^{\alpha}\delta _{j}^{i} \\ \noa
-\delta _{\beta }^{\alpha }\delta _{j}^{i} & 0
\end{array} \right)
$$
is a polysymplectic structure acting on $R^{np+np^{2}}$ with values in \
$R^{n+np}.$

Let $T_0 =[0,T^1]\times [0,T^2]\times \cdot \cdot \cdot  \times [0,T^p]$ and $dt = dt^1\cdot \cdot \cdot dt^p$ the volume element.
The action $\Psi$, whose Euler-Lagrange equations are the Hamilton
equations, is
$$
\Psi \left( u\right) = \int_{T_{0}} {\cal L} \left( t,u(t),\frac{\partial u}{\partial t}\right)dt, \quad \hbox{where} \quad u=\left( x,p\right),
$$
$$
{\cal L} \left( t,u(t),\frac{\partial u}{\partial t}\right) =-\frac{1}{2}
\left( \frac{\partial p_{i}^{\alpha }}{\partial t^{\alpha }}x^{i}-\frac{\partial x^{i}}{\partial t^{\alpha }}p_{i}^{\alpha }\right) -H\left(t,x,p\right) =
$$
$$
=-\frac{1}{2}\left( \frac{\partial p_{j}^{\alpha }}{\partial
t^{\alpha }},-\frac{\partial x^{j}}{\partial t^{\beta}}\right) \left(
\begin{array}{cc}
\delta _{ij} & o \\ \noa
o & \delta ^{\alpha \beta }\delta _{ij}
\end{array}
\right) \left(
\begin{array}{c}
x^{i} \\ \noa
p_{i}^{\alpha }
\end{array}
\right) -H\left( t,x\left( t\right), p\left( t\right) \right) =
$$
$$
=-\frac{1}{2} G\left( \left( \delta
\otimes J\right) \frac{\partial u}{\partial t},u\right) -H\left( t,u\right),
$$
where the scalar product is represented by the matrix $G=\left(
\begin{array}{cc}
\delta _{i j} & 0 \\ \noa
0 & \delta ^{\beta \alpha }\delta _{i j}
\end{array} \right) $
(standard Riemannian metric from $R^{n+np}$).
Indeed, the Euler-Lagrange equations produced by $\cal L$ are
$$
\frac{1}{2}\frac{\partial p_{i}^{\alpha }}
{\partial t^{\alpha }}=-\frac{1}{2}\frac{\partial p_{i}^{\alpha }}{\partial t^{\alpha }}-\frac{\partial H}{\partial x^{i}} \;
\hbox{,\,\,i.e.\,\,,}\; \frac{\partial p_{i}^{\alpha }}{\partial t^{\alpha }}=-\frac{\partial H}{\partial x^{i}}
$$
and
$$
\frac{\partial x^{i}}{\partial t^{\alpha }}=\frac{\partial H}{\partial p_{i}^{\alpha }}.
$$

It is well-known that the Legendre transform of a Lagrangian $L$ is the Hamiltonian $H$. Besides this classical duality,
we introduce a duality based on the Legendre transform of a Hamiltonian $H$.

The Legendre transformation $H^{\ast }\left( t,\cdot \right) $\ of the Hamiltonian $H\left(t,\cdot \right)$ is defined by the implicit formula
$$
H^{\ast }\left( t,v\right)
=G\left( u,v\right) -H\left( t,u\right),\,\, v=\nabla _{u}H\left(t,u\right),
$$
when $\nabla H\left(t,u\right)$ is invertible. Automatically,

$$ u=\nabla _{v}H^{\ast }\left( t,v\right),$$
so that $(\nabla H)^{-1}= \nabla H^{\ast}$.
If $H\left(t,\cdot \right)$ is of class $C^{1}$, strictly convex and has the property $\di\frac{H\left( t,u\right) }{\left| u\right| }\rightarrow \infty $
 for $\left|u\right| \rightarrow \infty $, then $H^{\ast}(t,\cdot)\in C^1(R^{n+np},R)$, according [2].
The norm $\left| u\right| $ is coming from the scalar product (Riemannian metric) $G$.

The Legendre transform $H^{\ast}(t,\cdot)$ of $H(t,\cdot)$ determines a new duality. Indeed, if we write $u = (x^i, p^{\alpha} _i)$,
the Hamiltonian equations can be written in the compact form (5). Setting $v = (y^i, q^{\alpha} _i)$ defined by $y^i =- x^i, q^{\alpha} _i = - p^{\alpha} _i$,
we obtain
$$
\left( \delta \otimes J\right) \di\frac{\partial v}{\partial t} =\left(
\begin{array}{c}
\di\frac{\partial q_{j}^{\alpha }}{\partial t^{\alpha }} \\ \noa
-\di\frac{\partial y^{j}}{\partial t^{\beta }}\end{array}\right)
= \left(
\begin{array}{c}
-\di\frac{\partial p_{j}^{\alpha }}{\partial t^{\alpha }} \\ \noa
\di\frac{\partial x^{j}}{\partial t^{\beta }}\end{array}\right)
 = -\left( \delta \otimes J\right) \di\frac{\partial u}{\partial t} = \nabla H(t,u).
$$
or equivalently

$$u = \nabla H^{\ast}(t,\left( \delta \otimes J\right) \di\frac{\partial v}{\partial t}).$$
if the Legendre transform $H^{\ast}(t,.)$ of $H(t,.)$ exists. On the other hand
$$
G\left( \left( \delta
\otimes J\right) \frac{\partial u}{\partial t},u\right) =\left( \frac{\partial p_{j}^{\alpha }}{\partial
t^{\alpha }},-\frac{\partial x^{j}}{\partial t^{\beta}}\right) \left(
\begin{array}{cc}
\delta _{ij} & o \\ \noa
o & \delta ^{\alpha \beta }\delta _{ij}
\end{array}
\right) \left(
\begin{array}{c}
x^{i} \\ \noa
p_{i}^{\alpha }
\end{array}
\right)
$$
$$
= \frac{\partial p_{i}^{\alpha }}{\partial t^{\alpha }}x^{i}-\frac{\partial x^{i}}{\partial t^{\alpha }}p_{i}^{\alpha }
$$
and consequently,
$$
G\left( \left( \delta
\otimes J\right) \frac{\partial v}{\partial t},v\right)=\frac{\partial q_{i}^{\alpha }}{\partial t^{\alpha }}y^{i}-\frac{\partial x^{i}}{\partial t^{\alpha }}q_{i}^{\alpha }=
G\left( \left( \delta
\otimes J\right) \frac{\partial u}{\partial t},u\right).
$$
Using the identity
$$
\di\frac{1}{2}G\left( -\left( \delta
\otimes J\right) \frac{\partial u}{\partial t},u\right)=\di\frac{1}{2}G\left( \left( \delta
\otimes J\right) \frac{\partial v}{\partial t},v\right) - G\left( \left( \delta
\otimes J\right) \frac{\partial v}{\partial t},v\right)
$$
$$
=\di\frac{1}{2}G\left( \left( \delta
\otimes J\right) \frac{\partial v}{\partial t},v\right) - G\left( \left( \delta
\otimes J\right) \frac{\partial u}{\partial t},u\right)
$$
$$
=\di\frac{1}{2}G\left( \left( \delta
\otimes J\right) \frac{\partial v}{\partial t},v\right) - G\left( \left( \delta
\otimes J\right) \frac{\partial v}{\partial t},u\right),
$$
the action
$$
\Psi \left( u\right) =\int_{T_{0}}[\frac{1}{2}G\left(
-\left( \delta \otimes J\right) \frac{\partial u}{\partial t},u\left(
t\right) \right) -H\left( t,u\left( t\right) \right)] dt.
$$
can be written as
$$
\Psi \left( u\right) =\int_{T_{0}}[\frac{1}{2}G\left(
\left( \delta \otimes J\right) \frac{\partial v}{\partial t},v\left(
t\right) \right) +G\left(\left( \delta \otimes J\right) \frac{\partial v}{\partial t} ,u\left( t\right) \right)-H\left( t,u\left( t\right) \right)] dt.
$$
In this way it appears the dual action
$$
\Phi \left( v\right)=\int_{T_{0}}[\frac{1}{2}G\left(
\left( \delta \otimes J\right) \frac{\partial v}{\partial t},v\left(
t\right) \right) + H^{\ast}\left( t,\left( \delta \otimes J\right) \frac{\partial v}{\partial t}\right)] dt. \eqno (6)
$$

\section{Extremals of dual action and multi-time Hamilton equations}

In this section we describe some connections between the periodical extremals of the
multi-time Hamilton dual action and the periodical solutions of the multi-time Hamilton equations.

The dual action, $\Phi \left( v\right)$,
is defined on a suitable space of T-periodic functions, where $T=(T^1,...,T^p)$. It uses the Lagrangian
$$
{\cal L}^{\ast }\left( t,v\left( t\right),
\frac{\partial v}{\partial t}\right) =\frac{1}{2}\left( \left( \delta
\otimes J\right) \frac{\partial v}{\partial t},v\left( t\right) \right)+H^{\ast}\left( t,\left( \delta \otimes J\right) \frac{\partial v}{\partial t}\right).
$$
Another useful property of the dual action is $\Phi \left( v+c\right) = \Phi \left(v\right)$, $\forall c\in R^{n+np}$. Indeed, for
$$
v^{\ast} = \left( \delta \otimes J\right) \frac{\partial v}{\partial t},\eqno (7)
$$
 we have
$$
\Phi \left( v+c\right) = \int_{T_{0}}[\left(
\left( \delta \otimes J\right) \frac{\partial \left( v+c\right) }{\partial t},v+c\right) +H^{\ast }\left( t,\left( v+c\right) ^{\ast }\right)] dt =
$$
$$
=\int_{T_{0}}[\left( \left( \delta \otimes
J\right) \frac{\partial v}{\partial t},v\right) +\left( \left( \delta
\otimes J\right) \frac{\partial v}{\partial t},c\right) +H^{\ast }\left(t,v^{\ast }\right) ]dt =
$$
$$
=\int_{T_{0}}[\left( \left( \delta \otimes
J\right) \frac{\partial v}{\partial t},v\right) +H^{\ast }\left( t,v^{\ast}\right)] dt=\Phi \left( v\right).
$$
Then we may restrict the search of the critical points of
$\Phi \left(v\right) $ to the set
$$
\widetilde{W}_{T}^{1,2}=\left\{ v\in
W_{T}^{1,2}\vert\int_{T_{0}}v\left( t\right) dt=0\right\}.
$$
There are situations when it is much easier to find the critical
points of the dual action $\Phi $ than those of the action $\Psi$.
This possibility was noticed by Clark in case of single-time
theory.

Between the periodical critical points of the dual action $\Phi $ and the periodical solutions of Hamilton multi-time equations,
there is a close connection, as we will show in the following result.

{\bf Theorem}. {\it Let $H:T_{0}\times R^{n+np}\rightarrow R$,
$\left(t,u\right) \rightarrow H\left( t,u\right)$ be a
Hamiltonian,  measurable in $t$ for any $u$ $\in R^{n+np}$ and
strictly convex and continuously differentiable in $u$ for any
$t\in T_{0}$. Suppose there are the constants $\alpha >0$, $\delta
>0$, and the functions $\beta$, $\gamma \in L^{2}\left(
T_{0,}R^{+}\right) $ so that
 for any $u\in R^{n+np}$ and any $t\in T_{0}$, to have
$$
\frac{\delta }{2}\left| u\right| ^{2} -\beta \left( t\right) \leq H\left(t,u\right) \leq\frac{\alpha}{2}\left| u\right| ^{2}+\gamma
\left( t\right). \eqno (8)
$$
Then, the dual action $\Phi $ is continuously differentiable on $\widetilde{W}_{T}^{^{1,2}}$. }

{\bf Proof.} From the hypothesis, the Hamiltonian $H\left( t,u\right) $ is strictly convex in $u$ and has
 the property $\di\frac{H\left( t,u\right) }{\left| u\right| }\rightarrow \infty $ when $\left| u\right| \rightarrow \infty $.
Then it follows (see [2]) that $H^{\ast }\left( t,.\right) \in C^{1}\left(R^{n+np},R\right)$. From the relation (8) and the equality
$$
\max_{u\in R^{n+np}}\left( \left( v^{\ast },u\right) -
\frac{\alpha }{2}\left| u^{2}\right| -\gamma \right) =\frac{\alpha ^{-1}}{2}\left| v^{\ast }\right| ^{2}-\gamma,\eqno (9)
$$
when $\alpha >0$, we obtain the inequalities
$$
\frac{\alpha ^{-1}}{2}\left| v^{\ast }\right| ^{2} -\gamma \left( t\right) \leq H^{\ast }\left( t,v^{\ast }\right) \leq
\frac{\delta^{-1}}{2}\left| v^{\ast }\right| ^{2} +\beta \left(t\right) \eqno (10)
$$
Taking into account the Proposition 2.2 from [2, page 34], with $N=n+np$ (and $n$, $p=q=2$), we obtain
$$
\left| \nabla H^{\ast }\left( t,v^{\ast }\right) \right| \leq \left( \left|v^{\ast }\right| +\beta \left( t\right) +
\gamma \left( t\right) \right)\delta ^{-1}+1\leq C_{1}\left| v^{\ast }\right| +C_{2}\left( \beta \left(t\right) +\gamma \left( t\right) +1\right), \eqno (11)
$$
where $C_{1,}\,\,C_{2}$ are positive constants. As $\beta +\gamma +1\in L^{2}(T_0,R^{+}$, from (10) and (11),
the Lagrangian
$${\cal L}^{\ast }\left( t,v(t),\di\frac{\partial v}{\partial t}\right) =\di\frac{1}{2}\left( \left( \delta \otimes J\right) \di\frac{\partial v}{\partial t},v\right) +
H^{\ast }\left( t,v^{\ast }\right)$$
 satisfies conditions similar to those of Theorem 1.4 from [2, page 10].
By consequence, the dual action $\Phi $ is continuously differentiable on $W_{T}^{1,2}$, and on $\widetilde{W}_{T}^{1,2}$.

{\bf Theorem}. {\it The same hypothesis as in the previous theorem. If $v\in \widetilde{W}_{T}^{^{1,2}}$\ is a critical
point for the action $\Phi $, periodical, with the period $T=\left(
T^{1},...,T^{p}\right)$, then the function $w\left( t\right) =\nabla
H^{\ast }\left( t,v^{\ast }\left( t\right) \right) $\ verifies the problem
$$
\left( \delta \otimes J\right) \frac{\partial w}{\partial t}+\nabla H^{\ast}\left( t,w\left( t\right) \right) =0,
\quad  w\left( 0\right) =w\left( T\right)
$$
and }
$$
v(t) =w(t)-\frac{1}{T^{1}...T^{p}}\int_{T_{0}}w\left(
t\right) dt.
$$

{\bf Proof}. We take $w=\left( z^{i},r_{i}^{\alpha }\right)$, $v^{\ast }=\left( v^{\ast i},v_{i}^{\ast \alpha }\right) $.
From the definition of the dual action, we have
$v=(-x^{i}, -p_{i}^{\alpha}),\;  \hbox{and} \;
v^{\ast }=\di\left( -\frac{\partial p_{i}^{\alpha }}{\partial t^{\alpha }},\frac{\partial x^{j}}{\partial t^{\beta}}\right).$
Then
$$
{\cal L}^{\ast }\left( t,v(t),\frac{\partial v}{\partial t}\right) =\frac{1}{2}\left( -\frac{\partial x^{i}}{\partial t^{\alpha }}p_{i}^{\alpha }+
\frac{\partial p_{i}^{\alpha }}{\partial t^{\alpha }}x^{i}\right) +H^{\ast}\left( t,-\frac{\partial p_{i}^{\alpha }}{\partial t^{\alpha }},
\frac{\partial x^{i}}{\partial t^{\alpha }}\right).
$$
If the function $v$ is a critical point for the action produced by ${\cal L}^{\ast }$, then the Euler-Lagrange
equations $\di\frac{\partial }{\partial t^{\alpha }}\di\frac{\partial
{\cal L}^{\ast }}{\partial x_{\alpha }^{i}}=\di\frac{\partial
{\cal L}^{\ast }}{\partial x^{i}}$, where $x_{\alpha }^{i}=\di\frac{\partial x^{i}}{\partial t^{\alpha }}$, are verified.
So, we obtain
$$
\frac{-1}{2}\frac{\partial p_{i}^{\alpha }}{\partial
t^{\alpha }}+\frac{\partial }{\partial t^{\alpha }}\frac{\partial H^{\ast }}{\partial v_{i}^{\ast \alpha }}=
\frac{1}{2}\frac{\partial p_{i}^{\alpha }}{\partial t^{\alpha }},$$
or
$$ \; \frac{\partial }{\partial t^{\alpha }}\frac{\partial H^{\ast }}
{\partial v_{i}^{\ast \alpha }}=\frac{\partial p_{i}^{\alpha }}{\partial t^{\alpha }}.
$$
Because ${}^t\left( z^{i},r_{i}^{\alpha }\right) =\nabla H^{\ast }\left(
t,v^{\ast i},v_{i}^{\ast \alpha }\right)$, we find $r_{i}^{\alpha }=\di\frac{\partial H^{\ast }}{\partial v_{i}^{\ast \alpha }}$ and then
$$
\frac{\partial r_{i}^{\alpha }}{\partial t^{\alpha }}=\frac{
\partial p_{i}^{\alpha }}{\partial t^{\alpha }}. \eqno (12)
$$
On the other side, the equality $w=\nabla H^{\ast }\left( t,-\di\frac{\partial p_{i}^{\alpha}}{\partial t^{\alpha }}, \di\frac{\partial x^{i}}{\partial t^{\alpha }}\right)$
and the duality produce
$$
{}^{t}\left( -\frac{\partial p_{i}^{\alpha }}{\partial
t^{\alpha }},\frac{\partial x^{i}}{\partial t^{\alpha }}\right) =\nabla H\left( t,w\right).
$$
By using the relation (12) we have
${}^t\left( -\di\frac{\partial r_{i}^{\alpha }}{\partial t^{\alpha }},\di\frac{\partial x^{i}}{\partial t^{\alpha }}\right) =\nabla H\left( t,w\right)$, and hence
$$
-\frac{\partial r_{i}^{\alpha }}{\partial t^{\alpha }}
=\nabla _{z^{i}}H\left( t,w\right). \eqno (13)
$$
From the definition of the function $H$ we find
$$
\frac{\partial z^{i}}{\partial t^{\alpha }}=
\frac{\partial H}{\partial r_{i}^{\alpha }}(t,z^i,r_i^{\alpha}). \eqno (14)
$$
From the relations (13) and (14) we obtain
$$
{}^t\left( -\di\frac{\partial
r_{i}^{\alpha }}{\partial t^{\alpha }}, \di\frac{\partial z^{i}}{\partial t^{\alpha }}\right) =
\nabla H\left( t,w\right) \; \hbox{or} \;  \left( \delta \otimes J\right) \frac{\partial w}{\partial t}=\nabla H\left(t,w\right).
$$
Consequently, $w$ is the solution of the equation $\left( \delta \otimes J\right) \di\frac{\partial w}{\partial t}+\nabla H^{\ast}\left( t,w\right) =0$.
 We consider $v$ a critical point for the dual action $\Phi(v) $ having mean value zero.
Then, also $v_{c}=v+c$ is a critical point for $\Phi $ because $\di\int_{T_{0}}\left( \left( \delta \otimes J\right) \di\frac{\partial v}{\partial t},c\right) dt=0$
(see the periodicity of $v$) and
$$
\Phi \left( v+c\right) =\int_{T_{0}}\left( \left( \delta
\otimes J\right) \frac{\partial \left( v+c\right) }{\partial t},v+c\right)+H^{\ast }\left( t,\left( v+c\right) ^{\ast }\right) dt=
$$
$$
=\int_{T_{0}}\left( \left( \delta \otimes J\right) \frac{
\partial v}{\partial t},v\right) +\left( \left( \delta \otimes J\right) \frac{\partial v}{\partial t},c\right) +H^{\ast }\left( t,v^{\ast }\right) dt=
$$
$$
=\int_{T_{0}}\left( \left( \delta \otimes J\right) \frac{
\partial v}{\partial t},v\right) +H^{\ast }\left( t,v^{\ast }\right) dt=\Phi \left( v\right).
$$
Because $v$ has the mean value zero, i.e., $\di\int_{T_{0}}v\left( t\right) dt=0$, we find
$$
\int_{T_{0}}\left( v_{c}\left( t\right)
-c\right) dt=0,
$$
and hence
$$
c=\frac{1}{T^{1}\cdot \cdot
\cdot T^{p}}\int_{T_{0}}u_{c}\left( t\right) dt.
$$
If $v_{c}$ is a critical point for $\Phi $, then $u_{c}=-v_c$ is a critical point for $\Psi $, i.e., it verifies the equation
$\left( \delta \otimes J\right)  \di\frac{\partial u_{c}}{\partial t}+\nabla H\left( t,u_{c}\right) =0$. By duality,
$$
u_{c}=\nabla H^{\ast }\left( t,{}^t\left( \delta \otimes J\right) \frac{
\partial u_{c}}{\partial t}\right) =\nabla H^{\ast }\left( t,-\left( \delta \otimes J\right) \frac{\partial u}{\partial t}\right) =\nabla H^{\ast}\left( t,v^{\ast }\right) =w.
$$
Then $v(t)=v_{c}(t)-c=
w(t)- \di\frac{1}{T^{1},...,T^{p}}
\di\int_{T_{0}}w\left( t\right) dt$. Because $v$ has the period $T=\left( T^{1},...,T^{p}\right) $ , we have
$v\left( 0\right) =v\left( T\right) $, so $w\left( 0\right) =w\left( T\right)$.

\bigskip

{\bf Acknowledgements}. The authors are grateful to Prof. Kostake
Teleman and Ana-Maria Teleman for their valuable comments on this
paper.
\bigskip

{\bf References}

\bigskip
[1] T. de Donder,\textit{Theorie invariative du calcul des variations}, 1935.

[2] J. Mawhin, M. Willem: \textit{Critical Point Theory and Hamiltonian Systems}, Springer-Verlag, 1989.

[3] M. Neagu, C. Udri\c{s}te: {\it From PDE Systems and Metrics to Geometric Multi-time Field Theories}, Seminarul de Mecanic\u a,
Sisteme Dinamice Diferen\c tiale, 79, Universitatea de Vest din Timi\c soara, 2001.

[4] P. J. Olver, The Equivalence Problem and Canonical Forms for Quadratic Lagrangians,
Advances in Applied Mathematics, 9 (1998), 226-257.

[5] A.-M. Teleman, C. Udri\c{s}te: \textit{On Hamiltonian Formalisms in Mathematical Physics}, 4-th International
Conference of Balkan Society of Geometers, Aristotle University of
Thessaloniki, 26-30 June, 2002.

[6] C.  Udri\c ste: \textit{Nonclassical Lagrangian Dynamics
and Potential Maps}, Proceedings of the Conference in
Mathematics in Honour of Professor Radu Ro\c sca on the
occasion of his Ninetieth Birthday, Katholieke University Brussel,
Katholieke University Leuven, Belgium, December 11-16, 1999; \\ http://xxx.lanl.gov/math.DS/0007060, (2000).

[7] C. Udri\c{s}te: \textit{Geometric Dynamics}, Kluwer Academic Publishers, Dodrecht/Boston /London, Mathematics and Its Applications, 513, 2000;
Southeast Asian Bulletin of Mathematics, Springer-Verlag 24 (2000).

[8] C. Udri\c ste: {\it Solutions of DEs and PDEs as
Potential Maps Using First Order Lagrangians}, Centenial Vr\^{a}nceanu,
Romanian Academy, University of Bucharest, June 30-July 4, (2000); http://xxx.lanl.gov/math.DS/0007061, (2000);
 Balkan Journal of Geometry and Its Applications 6, 1, 93-108, 2001.

[9] C. Udri\c{s}te: \textit{From integral manifolds and metrics to potential maps}, Atti del'Academia Peloritana dei Pericolanti,Classe 1 di Scienze Fis. Mat.
e Nat., 81-82, A 01006 (2003-2004), 1-14.

[10] C. Udri\c{s}te, M. Neagu: \textit{From PDE Systems and Metrics to Generalized Field Theories}, http://xxx.lanl.gov/abs/math.DG/0101207.

[11] C. Udri\c{s}te, M. Postolache: \textit{Atlas of Magnetic Geometric Dynamics}, Geometry Balkan Press, Bucharest, 2001.

[12] C. Udri\c{s}te, A.-M. Teleman: \textit{Hamilton Approaches of Fields Theory}, IJMMS, 57 (2004), 3045-3056;
 ICM Satelite Conference in Algebra and Related Topics,
University of Hong-Kong , 13-18.08.02.

University Politehnica of Bucharest \\
Department of Mathematics \\
Splaiul Independentei 313 \\
060042 Bucharest, Romania \\
email:udriste@mathem.pub.ro

\end{document}